\documentclass[10pt,a4paper,twoside]{amsart}
\usepackage{amsmath,bbm,amssymb,amsxtra}
\usepackage{enumerate}
\allowdisplaybreaks
\usepackage{epsfig}
\usepackage{ae}
\theoremstyle{definition}

\newtheorem{lemma}{Lemma}[section]
\newtheorem{proposition}[lemma]{Proposition}

\newtheorem{theorem}[lemma]{Theorem}

\newtheorem{definition}[lemma]{Definition}

\newcommand{\prop}[1]{\begin{proposition}\label{#1}
\sl }
\newcommand{\eprop}{\end{proposition}}
\newcommand{\thm}[1]{\begin{theorem}\label{#1}
\sl }
\newcommand{\ethm}{\end{theorem}}
\newcommand{\lem}[1]{\begin{lemma}\label{#1}
\sl }
\newcommand{\elem}{\end{lemma}}

\newcommand{\defin}[1]{\begin{definition}\label{#1}
\sl }
\newcommand{\edefin}{\end{definition}}

\newcommand{\beqno}{\begin{eqnarray*}}
\newcommand{\eeqno}{\end{eqnarray*}}
\newcommand{\beqla}[1] {\begin {eqnarray}\label{#1}}
\def\eeq {\end {eqnarray}}
\newcommand{\beq}{\begin {eqnarray}}

\newcommand{\real}{{\mathbb R}}

\newcommand{\integer}{{\mathbb Z}}

\newcommand{\complex}{{\mathbb C}}

\newcommand{\torus}{{\mathbb T}}

\newcommand{\D}{\mathbb{D}}
\newcommand{\E}{\mathbb{E}}

\newcommand{\C}{\mathbb{C}}

\newcommand{\Chat}{\widehat{\mathbb{C}}}

\newcommand{\prob}{{\mathbb P}}
\newcommand{\expec}{{{\mathbb E}\,}}

\newcommand{\diam}{{{\rm diam}\,}}

              \def\CO{{\mathcal O}}

\def\qq{ \begin{eqnarray} }
\def\qqq{ \end{eqnarray} }
\def\rr{ \begin{equation} }
\def\rrr{ \end{equation} }

\def\qq{ \begin{eqnarray} }
\def\qqq{ \end{eqnarray} }

\newcommand{\hf}{{_1\over^2}}

\newcommand{\BbbD}{\mathbb{D}}%

\newcommand{\KHI}{\raisebox{0.07cm}{$\chi$}}

\begin{document}
\title[Random Curves by Conformal Welding]
{Random Curves by Conformal Welding}

\author[K. Astala]{Kari Astala$^{1,4}$}
\address{University of Helsinki, Department of Mathematics and Statistics,
         P.O. Box 68 , FIN-00014 University of Helsinki, Finland}
\email{kari.astala@helsinki.fi}


\author[P. Jones]{Peter Jones$^{2}$}
\address{Department of Mathematics,
Yale University, 10 Hillhouse Ave, New Haven, CT, 06510, U.S.A.}
\email{jones@math.yale.edu}

\author[A. Kupiainen]{Antti Kupiainen$^{1,3}$}
\address{University of Helsinki, Department of Mathematics and Statistics,
         P.O. Box 68 , FIN-00014 University of Helsinki, Finland}
\email{antti.kupiainen@helsinki.fi}

\author[E. Saksman]{Eero Saksman$^1$}
\address{University of Helsinki, Department of Mathematics and Statistics,
         P.O. Box 68 , FIN-00014 University of Helsinki, Finland}
\email{eero.saksman@helsinki.fi}
\footnotetext[1]{Supported by the Academy of Finland, $^2$NSF,  $^3$ERC and $^4$the EU-network CODY. We thank M. Bauer, D. Bernard,  S. Rohde and S. Smirnov
for discussions}


\begin{abstract}
We construct a conformally invariant random family of closed curves in the plane by
 welding of random homeomorphisms of the unit circle given in terms of
 the exponential of Gaussian Free Field. We conjecture that
our curves are locally related to SLE$(\kappa)$ for $\kappa<4$.
\end{abstract}

\maketitle


\section{Introduction}
\label{se:intro}

A major breakthrough in the study of conformally
invariant random curves  in the plane occurred  when  O. Schramm \cite{Schra1}
introduced the Schramm-Loewner Evolution (SLE), a stochastic process
which describes such curves  growing in a fictitious time so that the curve of interest  is obtained as time tends to infinity. In this note we summarize a different construction \cite{us}
of random curves which is stationary i.e. the probability measure on curves is directly defined without
introducing an auxiliary time. We carry out this construction for closed curves, a case that is
not naturally covered by SLE.

Our construction is based on the idea of conformal welding which provides
a correspondence between Jordan curves on the extended plane 
$\Chat $ and a set of homeomorphisms of the
circle $\torus$. Given a Jordan curve  $\Gamma
\subset \Chat$, let
$$ f_+ : \D \to \Omega_+ \;\; \mbox{ and } \;\; f_- : \D_\infty \to \Omega_-
$$
be  a choice of Riemann mappings of the  unit disc $\BbbD$  and its
complement onto the components of  $\Chat \setminus \Gamma = \Omega_+ \cup \Omega_-$. By
Caratheodory's theorem $f_-$ and $f_+$ both extend continuously to
$\partial \D = \partial \D_\infty$, and thus \beqla{homeo11}
 \phi =  f_+ ^{-1} \circ f_-
\eeq is a homeomorphism of $\torus$. In the welding problem we are asked
to invert this process; given a homeomorphism $\phi: \torus \to \torus$ we
are to find a Jordan curve $\Gamma$ and conformal mappings $f_{\pm}$
onto the complementary domains $\Omega_{\pm}$ so that
(\ref{homeo11}) holds. It is clear that the welding problem, when solvable, has natural
conformal invariance attached to it; any image of the curve $\Gamma$
under a M\"obius transformation of $\Chat$ is equally a welding
curve. Similarly, if $\phi: \torus \to \torus$ admits a welding, then so do
all its compositions with M\"obius transformations of the disk.

We solve the welding problem for a   random, locally scale invariant set of
homeomorphisms $h_\omega: \torus\to \torus$, 
thereby obtaining a random set of Jordan curves.
  To define $h$, identify  the circle as $\torus =\real /\integer =[0,1)$.
Given a positive Borel measure $\tau$ without atoms we get a
  homeomorphism  $h:[0,1)\to [0,1)$ by
  \begin{equation}
h(t)=\tau ([0,t))/\tau ([0,1) )
\label{h1}
\end{equation}
 It was proposed  by the second
author some years ago that a natural class of homeomorphisms $h$ is obtained by
 taking $\tau$ formally proportional to $ e^{\beta X(t)}dt$ where
   $\beta\geq 0$ and  $X$
 is the Gaussian Free Field on the circle
i.e. the random field $X$  
 with covariance
\begin{equation}
\E\, X(t)X(t')=-\log |e^{2\pi it}-e^{2\pi it'}|.
\label{co1}
\end{equation}
For a rigorous definition one introduces a regularization $X_\varepsilon$ 
which is a.s. continuous if $\varepsilon>0$ and shows that almost surely the weak limit
of Borel measures
\begin{equation}\tau (dz)=\lim_{\varepsilon\to 0}
e^{\beta X_\varepsilon (z)}/\expec e^{\beta X_\varepsilon (z)} dz
\label{taudef}
\end{equation}
exists and defines a non-atomic random Borel measure on $[0,1]$.

Our main result is then: 

\thm{th:result} 
 For
$\beta^2< 2$ and almost surely in $\omega$, the formula (\ref{taudef})
defines a H\"older continuous 
circle homeomorphism, such that  the welding
problem has a  solution $\gamma$, where  $\gamma$ is a  Jordan curve  bounding a domain $\Omega = f_+(\BbbD)$ with
a  H\"older continuous Riemann mapping $f_+$. For a given $\omega$, the solution is unique up to a M\"obius map of the plane.
\ethm
The "critical (inverse)  temperature" $\beta_c=\sqrt{2}$ corresponds to loss
of continuity of the maps $h$. For $\beta\geq\beta_c$ the limit (\ref{taudef})
is zero almost surely. Based on theoretical physics \cite{bf} one may  conjecture that a corresponding limit of the normalized
measures $\tau_\varepsilon/\tau_\varepsilon([0,1])$ is non trivial also for
$\beta\geq\beta_c$ and atomic for $\beta>\beta_c$ thereby giving rise to
a  discontinuous map $h$. This phase transition is closely connected
to the one observed in two dimensional Liouville Quantum Gravity \cite{DuShe}
where a two dimensional version of our measure $\tau$ is considered. 

  We conjecture that the
 curves $\gamma$ locally "look like" SLE$(2\beta^2)$ (see also \cite{DuShe1} for
 arguments to this direction).  The case $\beta=\beta_c$,
  presumably corresponding to  SLE$(4)$, is not covered by our analysis.

 It would also be of interest to understand the connection of our weldings to those arising from stochastic flows studied in the interesting work \cite{AiMaTha}. In  \cite{AiMaTha} a program was set up for studying weldings that correspond to
   H\"older continuous homeomorphisms, but the boundary behaviour
 of the welding maps and hence the  existence and uniqueness  the welding was left open.
 
 \section{Beltrami equation}
\label{se:bel}

A powerful way to solve the welding problem goes by using
the Beltrami equation. Assume a homeomorphism $\phi:\torus\to\torus$ is
extended to a locally quasiconformal map $f:\D\to\D$, i.e.
 $f
\in C(\overline{\D})$ is a homeomorphism with $\nabla f$ locally integrable in $\D$ and satisfying
\beqla{homeo12} 
\frac{\partial f}{\partial \overline z} = \mu(z)
\frac{\partial f}{\partial  z}, \quad \quad \mbox{ for
} a.e.\;  z \in \D, 
\eeq 
with $\sup_{z\in K}|\mu(z)| < 1$ for $K \subset\subset D$.
One then considers the modified  equation 
\beqla{homeo14}
\frac{\partial F}{\partial \overline z} =
\KHI_{\D}(z) \, \mu(z) \; \frac{\partial
F}{\partial  z}, \quad \quad \mbox{ for } a.e. \; z \in \C. 
\eeq 

Suppose we can find a solution $F$ to (\ref{homeo14}) which
is a homeomorphism of $\Chat $. Then
$\Gamma = F(\torus)$ is a Jordan curve. Moreover, as
$\partial_{\overline z} F = 0$ for $|z| > 1$, we can set $ f_- :=
F|_{\D_\infty} \mbox{ and }   \Omega_- := F(\D_\infty)$ to define a
conformal mapping
$$ f_-: \D_\infty \to \Omega_-$$
To get the mapping $f_+$ note that both  $f$ and $F$ solve the Beltrami equation
(\ref{homeo12}) in the unit disk $\D$. 
By the uniqueness properties of equation (\ref{homeo12})
we conclude that
 \beqla{homeo31}
 F(z) = f_+ \circ f(z), \quad z \in \D,
\eeq for some conformal mapping $f_+ : \D = f(\D) \to \Omega_+ :=
F(\D)$. Then, on the unit circle,
\beqla{homeo32}
 \hskip40pt \phi(z)  = f|_{\torus}(z) =   f_+ ^{-1} \circ f_-(z), \quad \quad z \in \torus
\eeq
and we have found a solution to the welding problem.

To carry out this set of ideas we observe first that any homeomorphic self map $\phi$ of the circle can be extended to a locally quasiconformal map $f :\D\to\D$ via the Beurling-Ahlfors extension. However, a highly nontrivial problem remains: when does the auxiliary equation (\ref{homeo14}) have a locally quasiconformal
solution and when is this unique up to a conformal map?

A classical case where this question can be answered positively is the uniformly elliptic one where there is
an extension with $\|\mu\|_\infty<1$. This in turn will be true if 
 $\phi$  is
quasisymmetric. In our case these conditions {\it do not} hold,
and we are forced outside
the uniformly elliptic PDE's and need to study (\ref{homeo12}) with
strongly degenerate coefficients with only $|\mu(z)| < 1$ almost everywhere.

\section{Existence: Lehto method}
\label{se:lehto}

We use the method due to Lehto \cite{Le} to show   the existence of homeomorphic
solutions to (\ref{homeo14}). This approach is based on controlling the conformal moduli of
images of annular regions. 
To recall his result, define the distortion function 
$$K(z)  = \frac{1+|\mu(z) |}{1-|\mu(z) |}, 
$$
corresponding to the complex dilatation  $\mu=\mu (z)$.  Given an annulus 
$
A(w,r,R):=\{  z\in\complex : r<|z-w|<R\} 
$ define {\it the Lehto integral}:
\beqla{eq:4.7} L(w,r,R)
:=\int_r^R  \, \frac{1 }{
\int_0^{2\pi}
 {K}\left(w +\rho e^{i\theta}\right) d\theta} \; \frac{d\rho}{\rho}
\eeq
Lehto's theorem (see \cite[p.
584]{AsIwMa}) states then that if $\mu$ is  compactly supported with
 $|\mu(z)| < 1$ a.e.,  if $K(z)$ is locally integrable, and if
  for some $R_0>0$ the Lehto integral satisfies
 \begin{equation} \label{viimet27}  L(z,0,R_0) = 
 \infty,  \quad \mbox{for all}\quad z \in \C
 \end{equation}
then the Beltrami equation (\ref{homeo14})
admits a homeomorphic $W^{1,1}_{loc}$-solution $F:\C \to \C$.

We need actually a stronger result on the Lehto integrals to
obtain H\"older continuity of the solution. The Lehto integral
controls the geometric distortion
of an annulus under a locally quasiconformal map. Indeed,
given a bounded
(topological) annulus $A\subset\complex$ , with $E$ the bounded
component of $\complex\setminus A,$
 we denote by $D_O(A):=\diam (A)$ the outer diameter,
  and by $D_I(A):=\diam (E)$ the inner diameter of A.  It then
holds that for a quasiconformal map $f$
$$
D_I(f(A(w,r,R)))\leq 
{16}\exp\left(-2\pi^2 L(w,r,R)\right)D_O(f(A(w,r,R))).
$$
The equation (\ref{homeo14}) is solved by considering a regularized
uniformly elliptic equation where $\mu$ is replaced by $(1-\varepsilon)\mu$.
Since the corresponding solutions $F_{\varepsilon}$ are conformal in $\D_{\infty}$
 the outer radii $D_O(F_\varepsilon(A(w,r,1)))$ for $w\in\D$ are uniformly
 bounded by Koebe. Thus an estimate
  $$
  L(w,r,1)\geq a\log 1/r
  $$
  leads to $
D_I(F_\varepsilon(A(w,r,1)))\leq Cr^{2\pi^2 a}$ i.e. to H\"older continuity of
$F_\varepsilon$ uniformly in $\varepsilon$. Our main probabilistic estimate
then is
\thm{mainestimate}  Let $w\in\torus$ and let $\beta<\sqrt{2}$. Then there exists $b>0$ and $\delta_0 >0$  such that for
$\delta<\delta_0$ 
the Lehto integral satisfies the estimate
\begin{equation}\label{me}
\prob\bigl(L(w,2^{-N},1)<N\delta)\bigr)\leq
 2^{-(1+b )N}.
\end{equation}
\ethm
This estimate suffices to prove the existence and H\"older continuity of the solution
to  (\ref{homeo14}). First, only annuli centered at $w\in\partial \D$
need to be considered. Second, the $b>0$ allows us to cover, for each integer $N$,  
$\partial \D$ by balls $B_i$ of radii $2^{-(1+\hf b)N}$ such that for $\alpha>0$
$$
{\rm diam}(F(B_i))\leq C2^{-\alpha N}
$$
for all $i$ with probability $\CO(2^{-\hf bN})$. A Borel-Cantelli argument then
gives an a.s. H\"older continuity of $F$. 

\section{Uniqueness of the welding}
\label{se:uniqueness}

An important issue of the welding is its uniqueness, that the curve
$\Gamma$ is unique up to composing with a M\"obius transformation of
$\Chat$. This would follow from the uniqueness of solutions to the 
Beltrami equation (\ref{homeo14}), up to a M\"obius transformation. Unfortunately the control of the Lehto integrals given in Theorem \ref{mainestimate} alone is much too
weak to imply  this. 
However, in our case the uniqueness of solutions to the Beltrami
equation (\ref{homeo14})
 is equivalent to the conformal removability of the curve $F(\torus)$. 
 Indeed, suppose that we have two pairs 
 $f_\pm$ and $g_\pm$ of solutions to eq. (\ref{homeo11}). Then the formula
 $$\Psi(z) = \left\{\begin{array}{ll}
g_+ \circ  \left( f_+ \right)^{-1}(z) &\quad \mbox{if $z\in f_+(\D)$} \\
g_- \circ  \left( f_- \right)^{-1}(z)  &\quad \mbox{if $z\in
f_-(\D_\infty)$}
\end{array} \right.
 $$
 defines a homeomorphism of $\Chat$ that is conformal outside $\Gamma = f_\pm(\torus)$. 
 Since $\Gamma$ is a H\"older curve we can invoke the
  result of Jones and Smirnov in \cite{JoSmi} that H\"older curves are conformally removable
  i.e. that $\Psi$ extends conformally to the entire sphere. Thus it is a M\"obius transformation
  and uniqueness of the welding follows.

\section{A large deviation estimate}
\label{se:ld}

Theorem \ref{mainestimate}  follows from a large deviation estimate for
weakly correlated random variables. Let $\rho=2^{-p}$ where we choose
$p$ large. Let
$
L_k= L_{K_f}(w,\rho^{k},2\rho^{k})
$ so that $$L_{K_f}(w,2^{-Np},1)\geq \sum_{k=1}^NL_k.$$
For $p$ large $L_k$ are Lehto integrals in well separated
annuli (in logarithmic scale).
Estimate (\ref{me}) follows then from the inequality
\begin{equation}
{\prob}(  \sum_{k=1}^NL_k<N\delta)<\rho^{(1+b)N}.
 \label{prob}
\end{equation}
The bound (\ref{prob}) is a large deviation estimate and to prove it we establish two
facts: that (i) the random variables $L_k$ are (exponentially) weakly correlated
and (ii) uniformly in $k$, ${\prob}(L_k< \varepsilon)\leq C\varepsilon$.
These facts in turn rely on three ingredients: (a) an extension of $\phi$
to $f:\D\to \D$ with good local distortion bounds in terms of the random
measure $\tau$, (b) sharp probabilistic bounds for  $\tau$ and (c) a decomposition of the free field in 
terms of random fields localized in scale space.

For (a) we use the classical Beurling-Ahlfors extension \cite{BA}. We pave $\D$ by  
Whitney cubes $\{C_I\}_{I\in\mathcal D}$ indexed by dyadic intervals $I\subset\partial\D$
with ${\rm diam}(C_I)$ and ${\rm dist}(C_I, I)$ comparable to $|I|$. Then, 
extending some results by Reed on the Beurling-Ahlfors extension \cite{Re},  for $z\in C_I$
we have the  local distortion bound
\begin{equation}
K_f(z)\leq C \sum_{J,J'}{\tau(J)\over\tau(J')}
 \label{db}
\end{equation}
where 
 $J,J'$ run through dyadic intervals of size $ 2^{-4}|I|$ lying in $I$ and
 and its dyadic neighbours.
 The virtue of this bound is that the resulting lower bound for Lehto integral $L_k$ depends mostly
 on the ratios ${\tau(J)\over\tau(J')}$ for $J,J'$  dyadic intervals of size $\CO(2^{-kp})$
 and of distance $\CO(2^{-kp})$ from $w$. Thus we need to understand the sizes and
 mutual correlations of such ratios.
 
 For (b) we use results by Bacry and Muzy  \cite{BaMu} and Kahane \cite{Ka1} on multiplicative cascades (we refer the reader to \cite{BaMa} for an extensive discussion of  random multifractal measures). 
 The most crucial facts are that for $\beta <\sqrt{2}$ 
the measure $\tau$ is non-atomic and for any interval $I$ 
\begin{equation}
\tau(I)\in L^p(\omega), \ \ \ p\in(-\infty,2/\beta^2).
\label{lp1}
\end{equation}
Hence in particular the ratios in  (\ref{db}) are in $L^p(\omega)$ for $p\in[1,2/\beta^2)$.
 These facts are used in the proof of statement (ii) above. The fact that we may choose $p>1$
 is crucial for our analysis and is the source for the restriction to $\beta <\sqrt{2}$.
 
(c)  To understand the correlations between the $L_k$ i.e. between the ratios  ${\tau(J)\over\tau(J')}$
on scale $2^{-kp}$ we use a representation due to Bacry and Muzy  \cite{BaMu}  of
the free field $X$. It allows to decompose $X$ as 
 \begin{equation}
X=\sum_{k=0}^{\infty}\zeta_k
\label{decomp}
\end{equation}
where $\zeta_k$ are mutually independent a.s. continuous fields
with $\zeta_k(x)$ independent from $\zeta_k(y)$ for $|x-y|>\CO(2^{-kp})$.
This decomposition leads to the following  lower bound
\begin{equation}
L_n\geq m_n\exp({\sum_{k=0}^{n-1}2^{-{ap}(n-k)}t_{n,k}})(1+\sum_{k=n+1}^\infty 
2^{-{ap}(k-n)}\ell_{n,k})^{-1}
\label{Mnfinal}
\end{equation}
for $a>0$. The main contribution here is the scale $2^{-np}$ contribution $m_n$. The
positive random variables  $m_n$ are i.i.d.    
and satisfy the condition ${\prob}(m_n< \varepsilon)\leq C\varepsilon$.
Thus their sum $\sum m_n$ satisfies the estimate (\ref{prob}). 

The
corrections $t_{n,k}\geq 0$ and $\ell_{n,k}\geq 0$ represent correlations between
scale $2^{-np}$ and scale $2^{-kp}$ and are multiplied with exponentially small 
weights in
$|n-k|$. $t_{n,k}$ has gaussian tails:
$$
\prob (t_{n,k}> u )\leq
 ce^{-u^2/c}
 $$
 and the $\ell_{n,k}$ has a power law tail:
 $$
 \prob (\ell_{n,m}> \lambda)\leq C\lambda^{-q}
 $$
 for $q>1$.
Moreover, $t_{n,k}$ and $t_{n,'k'}$ are independent if $k\neq k'$
and $\ell_{n,m}$ and $\ell_{n',m'}$ are independent if $n>m'$ or $n'>m$.
These properties suffice to show that the estimate  (\ref{prob}) extends
from the $m_n$ to the $L_n$.

\end{document}